\documentclass{amsart}
\usepackage{amssymb}
\usepackage{amsfonts}
\usepackage{amsmath}
\usepackage{lastpage}
\usepackage[legalpaper,bookmarks=true,colorlinks=true,linkcolor=blue,citecolor=blue]{hyperref}
\usepackage{graphicx}%
\usepackage{fancyhdr}
\usepackage{color}
\usepackage[mathlines]{lineno}
\usepackage{lscape}
\usepackage{epsfig}
\usepackage{natbib}
\usepackage[]{listings}
\usepackage{geometry}
\usepackage{tgbonum}
\fontfamily{qcr}\selectfont
\newtheorem{theorem}{Theorem}
\theoremstyle{plain}

\numberwithin{equation}{section}

\newcommand{\Bin}{\bigskip \noindent}
\newcommand{\Bi}{\bigskip}
\newcommand{\Ni}{\noindent}

\begin{document}
\Large
\title{The Real-Valued Bochner integral and the Lebesgue-Like integration of real-valued measurable functions on $\mathbb{R}$}
\author{Gane Samb Lo}
\author{Lois Chinwendu Okereke}
\author{Fatima Doumbia}

\begin{abstract}
The Lebesgue-Like integral  of real-valued measurable functions (abbreviated as \textit{LL-integral}) is on one the most complete and appropriate integration Theory
 in Mathematics. Integrals are also defined in abstract spaces since Pettis (1938). In particular, Bochner integrals received much interest with very recent researches. It is very common to use the \textit{LL}integral while constructing other types of integrals, in particular the Bochner in Banach or in locally convex spaces. In this simple note, we prove that the Bochner integral and the LL-integral with respect to a finite measure $m$ are the same on $\mathbb{R}$. Applications of that equality may be useful in weak limits on Banach space.\\

\noindent $^{\dag\dag\dag}$ Gane Samb Lo.\\
LERSTAD, Gaston Berger University, Saint-Louis, S\'en\'egal (main affiliation).\newline
LSTA, Pierre and Marie Curie University, Paris VI, France.\newline
AUST - African University of Sciences and Technology, Abuja, Nigeria\\
gane-samb.lo@edu.ugb.sn, gslo@aust.edu.ng, ganesamblo@ganesamblo.net\\
Permanent address : 1178 Evanston Dr NW T3P 0J9,Calgary, Alberta, Canada.\\

\noindent $^{\dag\dag}$ Lois Chinwendu Okereke\\
AUST - African University of Sciences and Technology, Abuja, Nigeria\\

\noindent \noindent $^{\dag\dag\dag}$ Fatima Doumbia.\\
AUST - African University of Sciences and Technology, Abuja, Nigeria\\
\end{abstract}

\maketitle

\section{Introduction} \label{sec1}

\Ni This simple note focuses on the comparison between the LL-integral and the Bochner integral of a random variable $f$ defined on a measure space 
$(\Omega, \mathcal{A},m)$ and taking  values in $E=\mathbb{R}$, whenever they make sense. The LL-integral uses the order of $\mathbb{R}$ to completely describe the construction of the LL-integral, that we denote as

$$
\int_{(LL),\mathbb{R}} f \ dm \ \ or \ \ \int_{(LL)} f \ dm
$$

\Bi and call the (LL)-integral of $f$, by using a the three step method ($f$ elementary function, $f$ measurable and non-negative, $f$ simply measurable. The Bochner integral, that we denote as

$$
\int_{(Bo),\mathbb{R}} f \ dm \ \ or \ \ \int_{(Bo)} f \ dm 
$$

\Bin and call the (Bo)-integral of $f$,  which is constructed in a general Banach spaces ignores the order structure and use a two step methods.\\

\Ni Since both approaches are available on $\mathbb{R}$, we want to compare the two integral. So, we begin by summarizing the key elements of the construction
of the two types of integrals, restricting ourselves on the construction stages, to be in a position to make comparisons.\\

\subsection{Construction of the LL-integral} $ $\\

\noindent Assume that we have a measure space $(\Omega, \mathcal{A},m)$. We are going to construct the integral of a real-value measurable  function $f : (\Omega, \mathcal{A}) \rightarrow \overline{\mathbb{R}}$ with respect to the measure $m$ (that may take infinite values) denoted by

$$
\int_{(LL)} f \ dm=\int_{(LL),\Omega} f(\omega) \ dm(\omega) = \int_{(LL), \Omega} f(\omega) \text{ } m(d\omega).
$$

\noindent into three steps. The first step concerns non-negative functions among the class of elementary functions which have the general representation :

\begin{equation}
f=\sum_{1\leq i\leq k}\alpha _{i}\text{ }1_{A_{i}}, \ (\alpha_i \in \mathbb{R}_{+}, \ A_i \in \mathcal{A}, \ 1\leq i \leq k), \ k\geq 1,\label{elemFR}
\end{equation}

\noindent where the measurable sets $A_{i}$ are pairwise disjoint (\textit{pwd}). If it happens that the unions of the $A_i$, $1\leq i \leq p$ is not $\Omega$, we implicitly mean that $f=0$ on
the complement of $A_1+\cdots+A_p$.\\

\noindent The class of real-valued elementary functions is denoted by $\mathcal{E}(\omega, \mathcal{A}, \mathbb{R})$ and $\mathcal{E}^{+}(\omega,\mathcal{A}, \mathbb{R})$ stands for the subclass of non-negative functions of $\mathcal{E}(\omega,\mathcal{A}, \mathbb{R})$.\\

\noindent The expression of an elementary function, as expressed in Formula \eqref{elemFR} is not unique. But there exists one and only in which the coefficients $\alpha_i$ are disjoint, called the canonical representation. As a result, \textbf{that canonical representation is used with the summations of the $A_i$ covering $\Omega$} unless the contrary is specified.\\

\noindent Let us begin to describe the construction.\\

\bigskip  \noindent \textbf{Step 1M. Definition of the integral of a non-negative simple function : $f \in \mathcal{E}_{+}$}. \label{step1m}\\  

\noindent The integral of a non-negative simple function

\begin{equation}
f=\sum_{1\leq i\leq k}\alpha _{i}\text{ }1_{A_{i}}, \ (\alpha_i \in \mathbb{R}_{+}, \ A_i \in \mathcal{A}, \ 1\leq i \leq k), \ k\geq 1,\label{elemFR}
\end{equation}

\bigskip  \noindent is defined by

\begin{equation}
\int_{(LL)} f \ dm=\sum_{1\leq i\leq k} \alpha _{i} \ m(A_{i}). \ \ \label{IEF}
\end{equation}

\bigskip \noindent \textbf{Convention - Warning 1} In the definition \eqref{IEF}, the product $\alpha _{i} \ m(A_{i})$ is zero whenever $\alpha_i=0$, event if $m(A_i)=+\infty$.\\   

\noindent The definition \eqref{IEF} is coherent. This means that $\int_{(LL)} f \ dm$ does not depend on one particular expression of $f$.\\

\bigskip \noindent \textbf{step 2M. Definition of the integral for a non-negative measurable function}.\\ \label{step2m}

\noindent Let $f$ be any non-negative measurable function. By we have the following fact (see for example \textit{Point (03-23) in Doc 03-02 in Chapter 4 in \cite{ips-mestuto-ang})} : There exists a non-decreasing sequence $(f_n)_{n\geq 0}\subset \mathcal{E}_{+}$ such that

\begin{equation}
f_n \uparrow f \text{ as n } \uparrow +\infty. \label{approxPFEF}
\end{equation}

\bigskip \noindent  We define
 
\begin{equation}
\int_{(LL)} f \text{ }dm =\lim_{n\uparrow +\infty} \int_{(LL)} f_n  \ dm. \ \ \label{IPF}
\end{equation}

\bigskip \noindent  This definition \eqref{IPF} is also coherent since it does not depends of the sequence which is used in the definition. (See Chapter 4, \cite{ips-mestuto-ang}).\\

\bigskip \noindent \textbf{Step 3M. Definition of the integral for a measurable function}.\\ \label{step3m}

\noindent In the general case, the decomposition of $f$ into its positive part and its negative part is used as follows :

\begin{equation}
f=f^{+}-f^{-}, \ \  |f|=f^{+}+f^{-}, \ \ and \ \ f^{+}f^{-}=0, \ \label{PPPN}
\end{equation}

\bigskip \noindent where $f^{+}=max(0,f)$ and $f^{-}=max(0,-f)$, which are measurable, form the unique couple of functions such that Formulas \eqref{PPPN} holds.\\

\noindent By \textit{Step 2M}, the numbers 

$$
\int_{(LL)} f^+ \ dm \ \ and \int_{(LL)} f^- \ dm
$$

\bigskip \noindent exist in $\overline{\mathbb{R}}_{+}$. If one of them is finite, i.e.,

$$
\int_{(LL)} f^+ \ dm <+\infty \ \ \textbf{or} \int_{(LL)} f^- \ dm <+\infty,
$$

\bigskip \noindent we say that $f$ is quasi-integrable with respect to $m$ and we define

\begin{equation}
\int_{(LL)} f \ dm = \int f^+ \ dm \ - \ \int_{(LL)} f^- \ dm. \ \label{IGF}
\end{equation}

\bigskip \noindent \textbf{Warning} The integral of a real-valued and measurable function with respect to a measure $m$ exists \textit{if only if} : either it is of constant sign or the integral of its positive part or its negative part is finite.\\

\bigskip \noindent By extension, the \textbf{Integration of a mapping over a measurable set}. \noindent If $A$ is a measurable subset of $\Omega$ and $1_A f$ is quasi-integrable, we denote

$$
\int_{(LL),a} f \ dm = \int_{(LL)} 1_A f \ dm.
$$

\bigskip \noindent \textbf{Convention - Warning 2} The function $1_A f$ is defined by $1_A f(\omega)=f(\omega)$ of $\omega \in A$, and $zero$ otherwise.\\

\noindent The function $f$ is said to be (LL)-integrable if and only if the integral $\int f \ dm$ exist (in  $\mathbb{R}$) and is finite, i.e.,

$$
\int_{(LL)} f^+ \ dm <+\infty \ \ \textbf{and} \int_{(LL)} f^- \ dm <+\infty,
$$
  
\bigskip \noindent The set of all integrable functions with respect to $m$ is denoted

$$
\mathcal{L}^{1}(\Omega, \mathcal{A},m).
$$

\bigskip \noindent Let us move to the Bochner integral as given in \cite{mikusinski}.\\
 
\subsection{Bochner Integration on Banach Sapces} $ $\\

\noindent Given a finite measure space $\left(\Omega, \mathcal{A}, m \right)$, the Banach valued Bochner integral of measurable function

$$
f: \left(\Omega, \mathcal{A}, m \right) \longrightarrow (E, \mathcal{B}), \label{banachF}
$$ 

\bigskip \noindent where $(E,+, ., \left\|\circ\right\|_E)$ is a real Banach endowed with its Borel $\sigma$-algebra,  is defined through two steps. Below, the \textit{norm} of $f$, still denoted by $\left\|f\right\|_{E}$, is the measurable real-valued function defined by

$$
\forall \omega \in \Omega, \ \|f\|_{E}(\omega)=\|f(\omega)\|_{E}.
$$

\bigskip \noindent In the sequel, we denote by $\mathcal{E}(\Omega, \mathcal{A}, E)$ the class of all Banach-elementary functions.\\

\noindent \textbf{Step 1B}: The integral of an Banach-valued elementary function of the form

\begin{equation}
f = \sum_{j=1}^{p} x_{j} \ 1_{B_{j}}, \label{elFE} 
\end{equation}

\bigskip \noindent where $p\geq 1$, $x_{j} \in E$, $B_{j} \in \mathcal{A}$, $B_{1}+B_{2}+ \dots + B_{p} = \Omega$. At this stage,  only the linear structure of $E$ is used since $f(\omega)$ is a finite linear combination. Let us denote by 

\begin{equation}
\int_{(Bo)} f \ dm=\int_{(Bo} f \ dm =  \sum_{j=1}^{p} \ x_{j} m(A_{j}) \in E. \label{IBEF}
\end{equation}

\bigskip \noindent \textbf{Step 2B}: \label{step2B} A function \eqref{banachF} is Banach-Bochner integrable, denoted (Bo)-integrable, (denoted by $f \in \mathcal{L}^1(\Omega, \mathcal{A}, E, m)$) if and only if the two following conditions hold.\\

\noindent (a) There exists $ (f_n)_{n \geq 1} \subset \mathcal{E}(\Omega, \mathcal{A}, E)$ such that \\

\begin{equation}
\sum_{n \geq 1}\int_{((LL)} \left|f_n\right|_{E} \ dm < +\infty. \ \ \label{I1}
\end{equation}

\noindent (b) We have

\begin{equation}
f = \sum_{n \geq 1} f_n, \ m-a.e. \ \ \label{I2}
\end{equation}

\noindent If so, we write 

\begin{equation}
f \in S(f_n, n\geq1,\mathcal{E}(\Omega,\mathcal{A},E)). \label{attracSum}
\end{equation}

\bigskip \noindent and the Banach-valued Bochner integral is defined by 

\begin{equation}
\int_{(Bo)} f \ dm= \sum_{n\geq 1} \int_{(\Omega,E)} f_n \ dm. \label{integBochner}
\end{equation}

\bigskip \noindent A complete round on the subject is available in \cite{mikusinski}, where the consistency of the definitions of both steps has been been proved and the justification of both Formula \eqref{I1} and \eqref{I2}.\\

\noindent In the same paper, the space $\mathcal{L}^1(\Omega, \mathcal{A}, E, m)$ modulo the class of $m$-null-sets and denoted by $L^1(\Omega, \mathcal{A}, E, m)$
is proved to be a Banach space when endowed by the norm

$$
\left|f\right|_{L^1(\Omega, \mathcal{A}, E, m)}= \int_{(LL)} \left|f\right|_{E} \ dm.
$$

\bigskip \noindent Also, a dominated convergence theorem (DCT) is given therein.\\

\noindent As wll, the limits theory of sequences Bochner integrals is also very interesting. That theory is based on two key-ideas. First, we may replace $\mathcal{E}(\Omega,\mathcal{A},E)$ by $\mathcal{L}^1(\Omega,\mathcal{A},E)$ in Formula \eqref{attracSum} and still get the same decomposition of the integral as follows

\begin{theorem} Suppose that we have

$$
f \in S(f_n, n \geq 1,\mathcal{L}^{1}(\Omega,\mathcal{A},E)), \label{attracSumL1}
$$

\bigskip \noindent then $f$ is integrable and we have

\begin{equation}
\int_{(Bo)} f \ dm= \sum_{n\geq 1} int_{(\Omega,E)} f_n \ dm. \label{integBochner10}
\end{equation}
\end{theorem}

\noindent The second idea is that, in \textit{Step 1B}, for any $\eta>0$, we can choose a sequence $(f_n)_{n\geq 1}$ such that Formula \label{attracSum} holds and

\begin{equation}
 \left|f\right|_{E} \leq \sum_{n\geq 1}  \left|f_n\right|_{E} \leq \left|f\right|_{E}+\eta. \label{doubleApprox}
\end{equation}

\bigskip \noindent We have finished describing  two approaches. Let us proceed to their comparisons.\\
 
\section{Comparison of the two Integrals on $\mathbb{R}$}

\noindent In the previous section, we used the Real-valued Mapping \textit{(RVM)} integration scheme to get the Bochner integral in a complete normed space.\\

\noindent We already knew that the natural order on $\mathbb{R}$ was used in the general construction of the \textit{(RVM)} integration in the step \textit{2M} (see \pageref{step2m}) and we saw how that approach allowed to integrate with respect to an infinite measure.\\

\noindent Now, we are going to see that by restricting ourselves to finite measure, the Bochner and the \textit{(RVM)} integrals are exactly the same.\\

\noindent We suppose that the \textit{(RVM)} integration with respect to a finite measure $m$ is completely set. We have

\begin{theorem} \label{theoBochRVM} A real-valued measurable mapping and $m$-a.e. finite $f \ : \ (\Omega,\mathcal{A},m) \rightarrow \overline{\mathbb{R}}$ is 
(LL)-integrable if and only if its is (Bo)-integrable and its  its (LL)-integral and its (Bo)-integral coincide.
\end{theorem}

\bigskip \noindent \textbf{Proof of Theorem \ref{theoBochRVM}}. Let us consider a real-valued and measurable mapping and $m$-a.e. finite $f \ : \ (\Omega,\mathcal{A},m) \rightarrow \overline{\mathbb{R}}$. We have : for any elementary function, 
$$
\int_{(Bo)} f \ dm = \int_{(LL)} f \ dm.
$$

\Bin Let us proceed by step.\\

\Ni (a) Let $f$ be (LL)-integrable. Its positive part $f^+$ and negative part $f^-$ are integrable. In the context of $\mathbb{R}$, $f$ is limit of a non-decreasing sequence $(f_n^{(1)})_{n\geq 1}$ of non-negative elementary functions such that

$$
f_n^{(1)} \rightarrow f^+ \ as \ n \rightarrow +\infty \ \ and \ \ (\forall n\geq 1), \ |f_n^{(1)}|\leq f^+.
$$

\bigskip \noindent So, by the Dominated Convergence Theorem in the \textit{(RVM-MI)} scheme, we have

$$
\int_{(LL)} f_n^{(1)} \ dm \rightarrow \int_{(LL)} f^+ \ dm.
$$

\bigskip \noindent Now set $f_0=0$ and $h_n^{(1)}=f_n^{(1)}-f_{n-1}^{(1)}$ for $n\geq 1$. We have, for all $n\geq 1$,

$$
f_n^{(1)}=h_1^{(1)}+ \cdots + h_n^{(1)} 
$$

\bigskip \noindent It is clear that 

$$
\left( h_n^{(1)}\right)_{n\geq 1} \subset \mathcal{E}(\Omega,\mathcal{A},\mathbb{R})
$$

\bigskip \noindent and

$$
f^+=\sum_{1}^{+\infty} h_n^{(1)}
$$

\bigskip \noindent Further, we have for all $k\geq 1$,

\begin{eqnarray*}
\sum_{1}^{k} \int_{(LL)} |h_n^{(1)}| \ dm&=&\sum_{1}^{+\infty} \int_{(LL)} h_n^{(1)} \ dm\\
&=& \int_{(LL)}  f_k^{(1)} \ dm.
\end{eqnarray*}

\bigskip \noindent Hence, by taking the limit as $k\rightarrow +\infty$,

\begin{equation}
\sum_{1}^{+\infty} \int_{(LL)} |h_n^{(1)}| \ dm\leq  \int_{(LL)} f^+ \ dm<+\infty.
\end{equation}

\bigskip \noindent By doing the same for the negative part, we get a sequence

$$
\left( h_n^{(2)}\right)_{n\geq 1} \subset \mathcal{E}(\Omega,\mathcal{A},\mathbb{R})
$$

\bigskip \noindent such that 

$$
f^-=\sum_{1}^{+\infty} h_n^{(2)}
$$

\bigskip \noindent and

\begin{equation*}
\sum_{1}^{+\infty} \int_{(LL)} |h_n^{(2)}| \ dm\leq  \int_{(LL)} f^- \ dm<+\infty.
\end{equation*}

\bigskip \noindent Hence, by taking $h_n=h_n^{(1)}-h_n^{(2)}$, $n\geq 1$, we get

$$
\left( h_n\right)_{n\geq 1} \subset \mathcal{E}(\Omega,\mathcal{A},\mathbb{R}),
$$

$$
f=f^+-f^-=\sum_{1}^{+\infty} h_n
$$

\bigskip \noindent and

\begin{equation*}
\sum_{1}^{+\infty} \int_{(LL)} |h_n| \ dm\leq  \int_{(LL)} |f| \ dm<+\infty.
\end{equation*}

\bigskip
\noindent We conclude that

$$
f \in S(f_n, \ n\geq 1, \mathcal{E}(\Omega,\mathcal{A},\mathbb{R}))
$$

\bigskip
\noindent and hence $f$ is Bochner integrable and its Bochner integral is

$$
\int_{(Bo)} f \ dm=\sum_{n\geq 1} \int_{(LL)} h_n \ dm. 
$$

\bigskip \noindent Since

\begin{eqnarray*}
\sum_{n\geq 1} \int_{(LL)} h_n \ dm&=&\lim_{k\rightarrow +\infty} \sum_{1\leq n\leq k} \int_{(LL)} h_n \ dm\\
&=&\lim_{k\rightarrow +\infty}  \int_{(LL)} f_k^{(1)}-f_k^{(2)} \ dm\\
&=&\biggr(\lim_{k\rightarrow +\infty}  \int_{(LL)} f_k^{(1)}\biggr)-\biggr(\lim_{k\rightarrow +\infty}  \int_{(LL)} f_k^{(1)} \ dm\biggr)\\
&=& \int_{(LL)} f^+ \ dm -  \int_{(LL)} f^- \ dm\\
&=& \int_{(LL)} f \ dm,
\end{eqnarray*}

\noindent We have

$$
\int_{(Bo)} f \ dm = \int_{(LL)} f \ dm. \ \square
$$

\bigskip \noindent (a) Let $f$ be (Bo)-integrable with :

$$
f=\sum_{1}^{+\infty} f_n,
$$

$$
(f_n)_{n\geq 1} \subset \mathcal{E}(\Omega,\mathcal{A},\mathbb{R}),
$$

$$
\sum_{1}^{+\infty} \int_{(LL)} |f_n| \ dm <+\infty
$$

\bigskip \noindent and

$$
\int_{(Bo)} f \ dm =\lim_{k\rightarrow +\infty} \sum_{1\leq n\leq k} \int_{(LL)} f_n \ dm.
$$

\bigskip \noindent We take $g_k=\sum_{1\leq n\leq k} f_n \ dm$, $k\geq 1$. In the \textit{(LL)} scheme, the $g_k$'s are integrable (we recall that the measure is finite here!) and are bounded by 

$$
S=\sum_{n=1}^{+\infty} f_n, 
$$

\bigskip
\noindent which is (LL)-integrable. Then, by the dominated convergence theorem of the \textit{(LL)} scheme, we have that $f$ is integrable in that scheme and

$$
\int_{(LL)} f \ dm =\lim_{k\rightarrow +\infty} \sum_{1\leq n\leq k} \int_{(LL)} f_n \ dm=\int_{(Bo)} f \ dm .
$$ 

\bigskip \noindent The proof is over. $\blacksquare$\\

\noindent The main conclusions of that section are the following.\\

\noindent (1) The real-valued Bochner integral is exactly the modern real-valued integral with respect to a measure, provided the \textbf{measure is finite}.\\

\noindent (2) The Banach valued Bochner integral is an extension of the modern integral with respect to a finite measure to normed and complete space.\\

\noindent (3) The Bochner approach provides an alternative construction of the modern integral with respect to a finite measure, independently of the natural order of $\mathbb{R}$.

\end{document}